\documentclass[10pt]{amsart}
\usepackage{amsmath}
\usepackage{amssymb}
\usepackage{amsfonts}
\usepackage{amsthm}
\usepackage{amsxtra}
\usepackage{portland}
\usepackage{rotating}
\usepackage{nicefrac}
\usepackage{float}
\usepackage[all,web,arc,poly,dvips]{xy}
\usepackage{epic,eepic}
\usepackage{epsfig}
\usepackage[dvips]{color}
\usepackage{array}
\usepackage{pifont}
\usepackage{multirow}
\usepackage{graphics}
                         \theoremstyle{plain}

\newtheorem{corollary}{Corollary}[section]

\newtheorem{proposition}{Proposition}[section]
\theoremstyle{definition}

\theoremstyle{remark}

\newcommand{\ZZ}{\mathbb Z}
\newcommand{\NN}{\mathbb N}

\newcommand{\IIId}{${\rm III^{\prime}}\;$}

\begin{document}
                                                                                                    
\title[Boundary conditions associated with the Painlev\'e \IIId and V...]
{Boundary conditions associated with the Painlev\'e \IIId and V evaluations of 
some random matrix averages}
                                                                                                    
\author{P.J.~Forrester \and N.S.~Witte}
\address{Department of Mathematics and Statistics,
University of Melbourne,Victoria 3010, Australia}
\email{\tt p.forrester@ms.unimelb.edu.au} \email{\tt n.witte@ms.unimelb.edu.au}
                                                                                                    
\begin{abstract}
In a previous work a random matrix average for the Laguerre unitary ensemble, 
generalising the generating function for the probability that an interval
$ (0,s) $ at the hard edge contains $ k $ eigenvalues, was evaluated in terms 
of a Painlev\'e V transcendent in $ \sigma $-form. However the boundary conditions
for the corresponding differential equation were not specified for the full
parameter space. Here this task is accomplished in general, and the obtained
functional form is compared against the most general small $ s $ behaviour of the
Painlev\'e V equation in $ \sigma $-form known from the work of Jimbo. An analogous
study is carried out for the the hard edge scaling limit of the random matrix 
average, which we have previously evaluated in terms of a Painlev\'e \IIId
transcendent in $ \sigma $-form. An application of the latter result is given 
to the rapid evaluation of a Hankel determinant appearing in a recent work
of Conrey, Rubinstein and Snaith relating to the derivative of the Riemann
zeta function.
\end{abstract}
                                                                                                    
\maketitle

\section{Introduction}\label{S1}
\setcounter{equation}{0}
The Laguerre unitary ensemble ($ {\rm LUE}_N $) refers to the eigenvalue probability
density function (p.d.f.)
\begin{equation}
 \frac{1}{C_{N,a}}
 \prod^{N}_{l=1}\lambda_l^ae^{-\lambda_l}\prod_{1 \leq j<k \leq N}(\lambda_k-\lambda_j)^2 , 
\label{LUE_pdf}
\end{equation}
where
\begin{equation}
  C_{N,a} := \prod^{N-1}_{j=0}\Gamma(j+1)\Gamma(j+a+1) ,
\label{LUE_norm}
\end{equation}
with support on $ \lambda_l \in [0,\infty) $. For $ a \in \ZZ_{\geq 0} $, this
eigenvalue p.d.f. is realised by non-negative matrices $ X^{\dagger}X $ where 
$ X $ is an $ M\times N $ complex Gaussian matrix and $ a=M-N $. In our work
\cite{FW_2002a} the average over (\ref{LUE_pdf})
\begin{equation}
  \tilde{E}_{2,N}((0,s);a,\mu;\xi) := \frac{C_{N,a}}{C_{N,a+\mu}}
  \Big\langle \prod^N_{l=1}(1-\xi\chi^{(l)}_{(0,s)})(\lambda_l-s)^{\mu} \Big\rangle_{{\rm LUE}_N} ,
\label{LUE_E}
\end{equation}
where $ \chi^{(l)}_{J}=1 $ for $ \lambda_l \in J $, $ \chi^{(l)}_{J}=0 $ otherwise,
and the normalisation is chosen so that 
\begin{equation}
  \left.\tilde{E}_{2,N}((0,s);a,\mu;\xi)\right|_{s=0} := 1 ,
\label{LUE_Enorm}
\end{equation}
was characterised as a $ \tau $-function for the Painlev\'e V system. As a consequence,
it was shown that 
\begin{equation}
 W_N(s;a,\mu;\xi) 
  := s\frac{d}{ds}\log\left( s^{-N\mu}\tilde{E}_{2,N}((0,s);a,\mu;\xi) \right) ,
\label{LUE_sigma}
\end{equation}
satisfies the Jimbo-Miwa-Okamoto $ \sigma $-form of the Painlev\'e V equation
\begin{equation}
  (s\sigma''_V)^2 - \left[ \sigma_V - s\sigma'_V + 2(\sigma'_V)^2 +
  (\sum^3_{j=0}\nu_j) \sigma'_V \right]^2
  + 4 \prod^3_{j=0}(\nu_j + \sigma'_V) = 0 ,
\label{PV_sigma}
\end{equation}
with
\begin{equation}
  \nu_0 = 0, \quad \nu_1=-\mu, \quad \nu_2 = N+a, \quad \nu_3 = N, \quad
  \sum^3_{j=0}\nu_j = 2N+a-\mu .
\label{LUE_Vparam}
\end{equation}

For this to uniquely characterise $ W_N $, a boundary condition must be specified. 
However in \cite{FW_2002a} only in the cases $ \mu=0 $ and $ \mu=2 $ were a
boundary condition specified for general $ \xi $.

Also considered in \cite{FW_2002a} was the hard edge limiting average
\begin{equation}
  \tilde{E}^{\rm hard}_N(t;a,\mu;\xi) :=
  \lim_{N \to \infty} \Big(
  \frac{C_{N,a}}{C_{N,a+\mu}}\tilde{E}_{2,N}((0,\frac{t}{4N});a,\mu;\xi) \Big) .
\label{HE_Edefn}
\end{equation}
It was shown that
\begin{equation}
  \tilde{E}^{\rm hard}_2(t;a,\mu;\xi) = \exp\int_{0}^t u^{h}(s;a,\mu;\xi)\,\frac{ds}{s} ,
\label{HE_E}
\end{equation}
where setting
\begin{equation}
  u^{h}(s;a,\mu;\xi) =-(\sigma_{{\rm III}'}(s)+\frac{\mu(\mu+a)}{2}) ,
\label{HE_E.2}
\end{equation}
the function $ \sigma_{{\rm III}'}(s) $ satisfies the Jimbo-Miwa-Okamoto $ \sigma $-form
of the Painlev\'e \IIId equation
\begin{equation}
  (s\sigma''_{{\rm III}'})^2 - v_1v_2(\sigma'_{{\rm III}'})^2
  + \sigma'_{{\rm III}'}(4\sigma'_{{\rm III}'}-1)(\sigma_{{\rm III}'}-s\sigma'_{{\rm III}'})
  - \frac{1}{4^3}(v_1-v_2)^2 = 0 ,
\label{PIII_sigma}
\end{equation}
with parameters
\begin{equation}
  v_1 = a+\mu, \qquad v_2 = a-\mu .
\label{HE_IIIparam}
\end{equation}
Again only in the cases $ \mu=0 $ and $ \mu=2 $ were boundary conditions specified for 
general $ \xi $.

The aim of this work is to specify the boundary conditions relevant to both 
(\ref{LUE_sigma}) and (\ref{PV_sigma}) for general values of the parameters. In the 
case of (\ref{LUE_sigma}) this is done by writing the average (\ref{LUE_E}) in
its equivalent determinant form and evaluating the matrix elements in terms of
certain $ {}_1F_1 $ hypergeometric functions. With the small $ s $ behaviour 
of the matrix elements determined, it turns out that the determinant is such that
its corresponding small $ s $ behaviour can readily be deduced. The small $ s $
asymptotic form of (\ref{LUE_E}) then follows immediately. Scaling this 
asymptotic form as required by (\ref{HE_Edefn}) then gives the small $ t $
behaviour of $ \tilde{E}^{\rm hard}_N(t;a,\mu;\xi) $ and the small $ s $ behaviour
of $ u^{h}(s;a,\mu;\xi) $.

The general small $ s $ asymptotic form of the permitted solutions of 
(\ref{PV_sigma}) and (\ref{PIII_sigma}) have been given by Jimbo \cite{Ji_1982}.
As part of this study the boundary conditions found here are compared against
these general forms. It is found that in both cases only one of the two branches
permitted by the general solution is present in our random matrix problem.

As an application of our results we specify the rapid computation of the power
series expansion of a certain Hankel determinant of Bessel functions. The latter
is known form our work \cite{FW_2002a} to be a special case of 
$ \tilde{E}^{\rm hard}_N(t;a,\mu;\xi) $. The coefficients in the power series appear
in an asymptotic formula obtained recently by Conrey, Rubinstein and Snaith 
\cite{CRS_2005} for the integer moments of the derivative of the characteristic
polynomial of a unitary random matrix. This in turn has application to the study
of the derivative of the Riemann zeta function on the critical line.
 
\section{Small $ s $ Expansion of $ \tilde{E}_{2,N}((0,s);a,\mu;\xi) $}\label{S2}
\setcounter{equation}{0}
A standard result in random matrix theory, which in fact goes back to an identity
of Heine (see \cite{ops_Sz}) expresses the random matrix average (\ref{LUE_E}) as
a determinant
\begin{equation}
 \tilde{E}_{2,N}((0,s);a,\mu;\xi) 
  = \frac{N!C_{N,a}}{C_{N,a+\mu}}\det[w_{j+k}]_{j,k=0,\ldots,N-1} , 
\label{LUE_Hdet}
\end{equation}
where
\begin{equation}
  w_{n} 
  := \int^{\infty}_0 d\lambda\; (1-\xi\chi_{(0,s)})(\lambda-s)^{\mu}\lambda^{n+a}e^{-\lambda} .
\label{LUE_Hsym}
\end{equation}
Unless $ \mu $ is a non-negative integer (\ref{LUE_Hsym}) is not well defined for
$ s $ real and positive, which is the domain of interest. To remedy this, we note 
that simple manipulation gives
\begin{equation}
  w_{n} 
  := \int^{\infty}_s d\lambda\; (\lambda-s)^{\mu}\lambda^{n+a}e^{-\lambda}
     + (1-\xi)\int^{s}_0 d\lambda\; (\lambda-s)^{\mu}\lambda^{n+a}e^{-\lambda} ,
\label{LUE_Hsym.2}
\end{equation}
and in the second integral of this expression write 
$ (\lambda-s)^{\mu} = e^{\mu\log(\lambda-s)} $ with 
$ -\pi < \arg\log(\lambda-s) \leq \pi $. We then obtain
\begin{equation}
  w_{n} 
  := \int^{\infty}_s d\lambda\; (\lambda-s)^{\mu}\lambda^{n+a}e^{-\lambda}
     + (1-\xi)e^{\pi i\mu}\int^{s}_0 d\lambda\; (s-\lambda)^{\mu}\lambda^{n+a}e^{-\lambda} ,
\label{LUE_Hsym.3}
\end{equation}
which is well defined for $ \Re(\mu) > -1 $ and $ \Re(a) > -1 $ for 
$ s>0 $ with the additional constraint $ \Re(\mu+a) > -1 $ at $ s=0 $.

We seek the leading terms in the small $ s$ expansion of (\ref{LUE_Hsym.3}). 
These can be read off from an explicit evaluation in terms of the $ {}_1F_1 $
confluent hypergeometric function \cite{WW_1965}. 
\begin{proposition}
Subject to the conditions $ \Re(\mu) > -1 $, $ \Re(a) > -1 $,
$ \Re(\mu+a) > -1 $ and $ \mu+a\notin \ZZ_{\geq 0} $ we have
\begin{equation}
  w_{n} = a_n(s) + s^{n+\mu+a+1}b_n(s) ,
\label{LUE_Hsym.4}
\end{equation}
where $ a_n(s), b_n(s) $ are analytic about $ s=0 $ and given explicitly by
\begin{equation}
\begin{split}
  a_n(s) &=
  \Gamma(\mu+n+a+1)e^{-s}{}_1F_1(-a-n;-\mu-a-n;s) ,
  \\
  b_n(s) &=
  \frac{\Gamma(\mu+1)\Gamma(n+a+1)}{\Gamma(\mu+n+a+2)}
  \left( (1-\xi)e^{\pi i\mu}-\frac{\sin\pi a}{\sin\pi(\mu+a)} \right)
  \\
         & \qquad\times
  e^{-s}{}_1F_1(\mu+1;\mu+a+n+2;s) .
\end{split}
\label{Hsym_1F1}
\end{equation}
In particular, under the above conditions,
\begin{equation}
  w_{n} \mathop{\sim}\limits_{s \to 0} a_n(0)+sa'_n(0)+s^{n+\mu+a+1}b_n(0) ,
\label{LUE_Hsym.5}
\end{equation}
where
\begin{equation}
\begin{split}
  a_n(0)  &= \Gamma(\mu+n+a+1) ,
  \\
  a'_n(0) &= -\mu\Gamma(\mu+n+a) ,
  \\
  b_n(0) &=
  \frac{\Gamma(\mu+1)\Gamma(n+a+1)}{\Gamma(\mu+n+a+2)}
  \left( (1-\xi)e^{\pi i\mu}-\frac{\sin\pi a}{\sin\pi(\mu+a)} \right) .
\end{split}
\label{Hsym_1F1_exp}
\end{equation}
\end{proposition}
\begin{proof}
The results (\ref{LUE_Hsym.5}) and (\ref{Hsym_1F1_exp}) are immediate corollaries
of (\ref{LUE_Hsym.4}) and (\ref{Hsym_1F1}) and the fact that
\begin{equation*}
  {}_1F_1(\gamma;\alpha;s) = 1+\frac{\gamma}{\alpha}s+{\rm O}(s^2) .
\end{equation*}
To derive (\ref{LUE_Hsym.4}), we note that simple manipulation shows
\begin{equation*}
  \int^{\infty}_s d\lambda\; (\lambda-s)^{\mu}\lambda^{n+a}e^{-\lambda}
 = s^{a+n}e^{-s}\int^{\infty}_0 d\lambda\; (1+\lambda/s)^{n+a}\lambda^{\mu}e^{-\lambda} .
\end{equation*}
But with 
\begin{equation*}
  W_{k,m}(z) = \frac{z^ke^{-z/2}}{\Gamma(1/2-k+m)}
   \int^{\infty}_0 dt\; (1+t/z)^{k-1/2+m}t^{-k-1/2+m}e^{-t} ,
\end{equation*}
specifying the Whittaker function, it is known that \cite{WW_1965}
\begin{equation*}
  W_{k,m}(z)
  = \frac{\Gamma(-2m)}{\Gamma(1/2-k-m)}M_{k,m}(z)+\frac{\Gamma(2m)}{\Gamma(1/2-k+m)}M_{k,-m}(z)
\end{equation*}
where
\begin{equation*}
  M_{k,m}(z) = z^{m+1/2}e^{-z/2}{}_1F_1(1/2-k+m;2m+1;z) .
\end{equation*}
Consequently
\begin{multline}
  \int^{\infty}_s d\lambda\; (\lambda-s)^{\mu}\lambda^{n+a}e^{-\lambda}
 =  \Gamma(\mu+a+n+1)e^{-s}{}_1F_1(-a-n;-\mu-a-n;s) 
 \\
  + \frac{\Gamma(\mu+1)\Gamma(-\mu-a-n-1)}{\Gamma(-a-n)}
    s^{\mu+a+n+1}e^{-s}{}_1F_1(\mu+1;\mu+a+n+2;s) .
\label{ws6}
\end{multline}
The left-hand side of (\ref{ws6}) exists for $ \Re(\mu)>-1 $ if $ s>0 $ and
$ \Re(\mu+a)>-1 $ if $ s=0 $, whereas the right-hand side is valid in 
this parameter domain except for $ \mu+a+n\in\ZZ_{\geq 0} $, and in this case the individual
terms have a simple pole at $ a+n\notin\ZZ_{\geq 0} $ or are undefined when $ a+n\in\ZZ_{\geq 0} $. 
Needless to say the sum of the terms on the right-hand side has the same analytic
character as the left-hand side.

Regarding the second integral in (\ref{LUE_Hsym.3}), we first note that a simple 
change of variables gives
\begin{equation*}
  \int^{s}_0 d\lambda\; (s-\lambda)^{\mu}\lambda^{n+a}e^{-\lambda}
  = s^{n+1+a+\mu}e^{-s}\int^{1}_0 dx\; (1-x)^{n+a}x^{\mu}e^{sx} .
\end{equation*}
But the integral on the right hand side is the Euler integral representation of
the $ {}_1F_1 $ function, which shows
\begin{multline}
 \int^{s}_0 d\lambda\; (s-\lambda)^{\mu}\lambda^{n+a}e^{-\lambda} 
 \\
  = \frac{\Gamma(\mu+1)\Gamma(a+n+1)}{\Gamma(\mu+a+n+2)}
    s^{\mu+a+n+1}e^{-s}{}_1F_1(\mu+1;\mu+a+n+2;s) .
\label{ws7}
\end{multline}
This latter relation is valid for $ \Re(\mu)>-1 $ and $ \Re(a)>-1 $ when
$ s>0 $.
Substituting (\ref{ws6}) and (\ref{ws7}) in (\ref{LUE_Hsym.3}) and using the 
appropriate gamma function identities gives (\ref{LUE_Hsym.4}), (\ref{Hsym_1F1}).
\end{proof}

When $ \mu+a \in\ZZ_{\geq 0} $ we have to consider two exceptional cases where one of the 
hypergeometric functions are not defined - the first when $ a+n \in\ZZ_{\geq 0} $ 
for which the hypergeometric function is indeterminate, and the second when
$ a+n \notin\ZZ_{\geq 0} $ and the hypergeometric function has a simple pole.
These two cases can be recovered by taking suitable limits and we just state the
final results.

\begin{proposition}
When $ \mu+a=j \in\ZZ_{\geq 0} $ and $ a+n=k \in\ZZ_{\geq 0} $ with $ n+j\geq k $ 
we have
\begin{multline}
 w_n = k!e^{-s}\Bigg\{  \sum^k_{l=0}\frac{(n+j-l)!}{(k-l)!l!}s^l
 \\
  + (-1)^{n+j+k}(1-\xi)\frac{(n+j-k)!}{(n+j+1)!}s^{n+j+1}
    {}_1F_1(n+j+1-k;n+j+2;s) \Bigg\} ,
\label{Hsym_indeterm}
\end{multline}
and to leading order in small $ s $ we have
\begin{equation}
  w_{n} \mathop{\sim}\limits_{s \to 0} 
  (n+j)! - (n+j-k)(n+j-1)!s + (-1)^{n+j+k}(1-\xi)\frac{(n+j-k)!k!}{(n+j+1)!}s^{n+j+1} .
\label{Exp_indeterm}
\end{equation}
\end{proposition}
Note that the condition $ n+j\geq k $ is the same as $ \mu \geq 0 $, which falls
within the domain of interest. The key difference of (\ref{Exp_indeterm}) with
(\ref{LUE_Hsym.5}) and (\ref{Hsym_1F1_exp}) is that the non-analytic term is
now polynomial and the second part of this term is absent having been cancelled by
a counterbalancing term.

\begin{proposition}
When $ \mu+a=j \in\ZZ_{\geq 0} $ and $ a+n \notin\ZZ_{\geq 0} $ we have 
\begin{multline}
 w_n = e^{-s}\Bigg\{ \sum^{n+j}_{l=0}\frac{(-a-n)_l(n+j-l)!}{l!}(-s)^l 
 \\
  +\frac{\Gamma(\mu+1)\Gamma(a+n+1)}{(n+j+1)!}(1-\xi)e^{i\pi\mu}s^{n+j+1}
   {}_1F_1(\mu+1;n+j+2;s)
 \\
  + (-1)^j\frac{\sin\pi a}{\pi}\frac{\Gamma(\mu+1)\Gamma(a+n+1)}{(n+j+1)!}s^{n+j+1}
 \\
 \times
    \sum^{\infty}_{l=0}\left[\psi(l+1)+\psi(n+j+l+2)-\psi(\mu+l+1)-\log s \right]
          \frac{(\mu+1)_l}{(n+j+2)_l}\frac{s^l}{l!} \Bigg\} ,
\label{Hsym_pole}
\end{multline}
and its leading order behaviour for small $ s $ is
\begin{multline}
  w_{n} \mathop{\sim}\limits_{s \to 0} 
  (n+j)!+(a-j)(n+j-1)!s 
  \\
  +\frac{(a-j)_{n+j+1}}{(n+j+1)!}s^{n+j+1}
   \left[\frac{\pi e^{-i\pi a}}{\sin\pi a}(1-\xi)+\psi(1)+\psi(n+j+2)-\psi(\mu+1)-\log s \right] .
\label{Exp_pole}
\end{multline}
\end{proposition}
The expansion (\ref{Exp_pole}) differs significantly from (\ref{LUE_Hsym.5}) and 
(\ref{Hsym_1F1_exp}) because of the presence of logarithmic terms which now 
replace the non-analytic contributions of the generic case.

\begin{corollary}
Under generic conditions on $ \mu+a $ we have
\begin{multline}
  \det[w_{j+k}]_{j,k=0,\ldots,N-1} 
  \\
  = \det[\Gamma(\mu+a+1+j+k)]_{j,k=0,\ldots,N-1}
  \\
    - \mu s\det[\Gamma(\mu+a+j) \;\; \Gamma(\mu+a+1+j+k)]_{{j=0,\ldots,N-1}\atop
                                                         {k=1,\ldots,N-1}}
    + {\rm O}(s^2)
  \\
    + s^{\mu+a+1}b_0(0)\det[\Gamma(\mu+a+3+j+k)]_{j,k=0,\ldots,N-2}
      \left\{1+{\rm O}(s)\right\}
  \\
    + {\rm O}(s^{2(\mu+a+1)}) .
\label{Hdet_exp}
\end{multline}
\end{corollary}
\begin{proof}
According to (\ref{LUE_Hsym.5})
\begin{multline*}
  \det[w_{j+k}]_{j,k=0,\ldots,N-1} 
  \\
  \mathop{\sim}\limits_{s \to 0}
  \det[a_{j+k}(0)+sa'_{j+k}(0)+s^{\mu+a+1+j+k}b_{j+k}(0)]_{j,k=0,\ldots,N-1}
  \\
  \mathop{\sim}\limits_{s \to 0}
  \det[a_{j+k}(0)]_{j,k=0,\ldots,N-1}+s[s]\det[a_{j+k}(0)+sa'_{j+k}(0)]_{j,k=0,\ldots,N-1}
  \\
  + s^{\mu+a+1}b_0(0)\det[a_{j+k+2}(0)]_{j,k=0,\ldots,N-2} ,
\end{multline*}
where $ [s]P(s) $ denotes the coefficient of $ s $ in $ P(s) $. Recalling the 
explicit formula for $ a_n(0) $ as given in (\ref{Hsym_1F1_exp}) we obtain the 
constant term and the term proportional to $ s^{\mu+a+1} $ in (\ref{Hdet_exp}).
It remains to compute the coefficient of $ s $, which according to 
(\ref{Hsym_1F1_exp}) has the explicit form
\begin{equation}
   [s]\det[\Gamma(\mu+a+1+j+k)-\mu s\Gamma(\mu+a+j+k)]_{j,k=0,\ldots,N-1} .
\label{ws9}
\end{equation}
Using the linearity formula
\begin{equation*}
   \det[{\bf a}_1 \cdots {\bf a}_j+{\bf b}_j \cdots {\bf a}_n]
   = \det[{\bf a}_1 \cdots {\bf a}_j \cdots {\bf a}_n]
    +\det[{\bf a}_1 \cdots {\bf b}_j \cdots {\bf a}_n] ,
\end{equation*}
where the $ \bf{a} $'s and $ \bf{b} $'s are column vectors, on each column of
the determinant we see that of the terms proportional to $ s $ only the one
obtained from expanding the first column in non-zero (all the rest result in 
two identical columns), and the determinant given by (\ref{Hdet_exp}) results.
\end{proof}

It remains to evaluate the determinants. For this task we make use of the identity
\cite{Nd_2004}
\begin{equation*}
   \det[\Gamma(z_k+j)]_{j,k=0,\ldots,n-1}
   = \prod^{n-1}_{j=0}\Gamma(z_j)\prod_{0\leq j<k\leq n-1}(z_k-z_j) .
\end{equation*}
After straightforward manipulations, gamma function evaluations of all the 
determinants in (\ref{Hdet_exp}) can be obtained. Substituting in (\ref{LUE_Hdet}),
and recalling that the normalisation is such that at $ s=0 $ $ \tilde{E}_N $ 
is equal to unity, we obtain the sought small $ s $ expansion of $ \tilde{E}_N $
and thus $ W_N $ valid for general values of the parameters.

\begin{proposition}\label{LUE_exp}
For $ \Re(\mu)>-1 $, $ \Re(a)>-1 $ and $ \mu+a\notin\ZZ_{\geq 0} $ we have 
\begin{multline}
  \tilde{E}_{2,N}((0,s);a,\mu;\xi) = 1-\frac{\mu N}{\mu+a}s+{\rm O}(s^2)
  \\
  +\frac{\Gamma(\mu+1)\Gamma(a+1)\Gamma(\mu+a+N+1)}{\Gamma^2(\mu+a+2)\Gamma(\mu+a+1)\Gamma(N)}
   \left( (1-\xi)e^{\pi i\mu}-\frac{\sin\pi a}{\sin\pi(\mu+a)} \right) 
   s^{\mu+a+1}\left\{1+{\rm O}(s)\right\}
  \\ +{\rm O}(s^{2(\mu+a+1)}) ,
\label{LUE_Eexp}
\end{multline}
and consequently
\begin{multline}
  W_{N}(s;a,\mu;\xi) = -N\mu-\frac{\mu N}{\mu+a}s+{\rm O}(s^2)
  \\
  +\frac{\Gamma(\mu+1)\Gamma(a+1)\Gamma(\mu+a+N+1)}{\Gamma(\mu+a+2)\Gamma^2(\mu+a+1)\Gamma(N)}
   \left( (1-\xi)e^{\pi i\mu}-\frac{\sin\pi a}{\sin\pi(\mu+a)} \right) 
   s^{\mu+a+1}\left\{1+{\rm O}(s)\right\}
  \\ +{\rm O}(s^{2(\mu+a+1)}) .
\label{LUE_Wexp}
\end{multline}
\end{proposition}
In the first exceptional case $ \mu+a=j \in\ZZ_{\geq 0} $ and $ a=k \in\ZZ_{\geq 0} $ 
with $ j\geq k $ one can still use (\ref{LUE_Eexp}) but omitting the term involving
the ratio of sines, in the case $ j=0 $, or the whole term if $ j>0 $.
The situation of the other exceptional case $ \mu+a=j \in\ZZ_{\geq 0} $ and 
$ a \notin\ZZ_{\geq 0} $ is more complicated and more so for larger $ j $, and we 
only give the examples of $ j=0,1 $.
\begin{proposition}
For $ \Re(\mu)>-1 $, $ \Re(a)>-1 $ with $ \mu+a=0 $ we have 
\begin{multline}
\tilde{E}_{2,N}((0,s);a,\mu=-a;\xi) = 1
  \\
  +\Big\{-1 + \frac{\pi a}{\sin\pi a}e^{-i\pi a}(1-\xi)
   +a\left[2\psi(2)+\psi(1)-\psi(1-a)-\psi(N+1)-\log s \right] \Big\}Ns
  \\ + {\rm o}(s) .
\label{LUE_Eexp_pole0}
\end{multline}
For $ \mu+a=1 $ we have
\begin{multline}
\tilde{E}_{2,N}((0,s);a,\mu=1-a;\xi) = 1+(a-1)Ns
  \\
  +\frac{a(a-1)}{4}\Big\{\frac{\pi}{\sin\pi a}e^{-i\pi a}(1-\xi)
   +2\psi(3)+\psi(2)-\psi(2-a)-\psi(N+2)-\log s \Big\}(N+1)Ns^2
  \\ + {\rm o}(s^2) .
\label{LUE_Eexp_pole1}
\end{multline}
\end{proposition}

\section{Comparison with the Jimbo solution}\label{S3}
\setcounter{equation}{0}
The small $ s$ expansion of the most general solution permitted by (\ref{PV_sigma}),
or more precisely its corresponding $ \tau $-function (see (\ref{PV_tau}) below) has 
been determined by Jimbo \cite{Ji_1982}. However in \cite{Ji_1982} the equation
(\ref{PV_sigma}) is not treated directly. Instead the discussion is based on
the equation
\begin{multline}
  (s\zeta'')^2 - \left[ \zeta - s\zeta' + 2(\zeta')^2 -
  (2\theta_0+\theta_{\infty}) \zeta' \right]^2 \\
  + 4\zeta'(\zeta'-\theta_0)(\zeta'-\frac{1}{2}(\theta_0-\theta_s+\theta_{\infty}))
    (\zeta'-\frac{1}{2}(\theta_0+\theta_s+\theta_{\infty}))
  = 0 ,
\label{PV_zeta}
\end{multline}
and the small $ s $ behaviour of the corresponding $ \tau $-function $ \tau_V(s) $,
specified by the the requirement that
\begin{equation}
  \zeta(s) = s\frac{d}{ds}\log\tau_V(s)+\frac{1}{2}(\theta_0+\theta_{\infty})s
          +\frac{1}{4}[(\theta_0+\theta_{\infty})^2-\theta_s^2] , 
\label{PV_tau}
\end{equation}
was determined.

Comparison of (\ref{PV_zeta}), (\ref{PV_tau}) with (\ref{PV_sigma}), 
(\ref{LUE_sigma}) shows that for the parameters (\ref{LUE_Vparam})
\begin{equation}
  \tilde{E}_{2,N}((0,s);a,\mu;\xi) = s^{N^2+N(a+\mu)}e^{-(N+a/2)s}\tau_V(s) ,
\label{LUE_Vtau}
\end{equation}
while in general
\begin{equation}
  \theta_0 = -\nu_1, \quad \theta_s = \nu_2-\nu_3, 
  \quad \theta_{\infty} = \nu_1-\nu_2-\nu_3 .
\label{Vnu_theta}
\end{equation}
Note that for the parameters (\ref{LUE_Vparam}) we thus thus have
\begin{equation}
  \theta_0 = \mu, \quad \theta_s = a, \quad \theta_{\infty} = -2N-a-\mu .
\label{LUE_theta}
\end{equation}

The relevant result from \cite{Ji_1982} can now be presented. It states that the
most general small $ s $ behaviour of $ \tau_V(s) $ permitted by the equation
(\ref{PV_zeta}) is
\begin{multline}
  \tau_V(s) = Cs^{(\sigma^2-\theta^2_{\infty})/4}
       \Bigg\{ 1
               - \frac{\theta_{\infty}(\theta^2_s-\theta^2_0+\sigma^2)}
                     {4\sigma^2}s \\
   + u
      \frac{\Gamma^2(-\sigma)}{\Gamma^2(2+\sigma)}
      \frac{\Gamma(1+\frac{\displaystyle\theta_s+\theta_0+\sigma}{\displaystyle 2})
            \Gamma(1+\frac{\displaystyle\theta_s-\theta_0+\sigma}{\displaystyle 2})
            \Gamma(1+\frac{\displaystyle\theta_{\infty}+\sigma}{\displaystyle 2})}
           {\Gamma(\frac{\displaystyle\theta_s+\theta_0-\sigma}{\displaystyle 2})
            \Gamma(\frac{\displaystyle\theta_s-\theta_0-\sigma}{\displaystyle 2})
            \Gamma(\frac{\displaystyle\theta_{\infty}-\sigma}{\displaystyle 2})}
      s^{1+\sigma}
   \\
   + \frac{1}{u}
      \frac{\Gamma^2(\sigma)}{\Gamma^2(2-\sigma)}
      \frac{\Gamma(1+\frac{\displaystyle\theta_s+\theta_0-\sigma}{\displaystyle 2})
            \Gamma(1+\frac{\displaystyle\theta_s-\theta_0-\sigma}{\displaystyle 2})
            \Gamma(1+\frac{\displaystyle\theta_{\infty}-\sigma}{\displaystyle 2})}
           {\Gamma(\frac{\displaystyle\theta_s+\theta_0+\sigma}{\displaystyle 2})
            \Gamma(\frac{\displaystyle\theta_s-\theta_0+\sigma}{\displaystyle 2})
            \Gamma(\frac{\displaystyle\theta_{\infty}+\sigma}{\displaystyle 2})}
      s^{1-\sigma}
   \\
   + {\rm O}(|s|^{2(1-\Re(\sigma))}) \Bigg\} ,
\label{PV_sExp_0}
\end{multline}
where $ C $ is a normalisation constant, while $ u $ and $ \sigma $ are arbitrary
parameters. The above result was derived subject to the conditions 
$ \theta_0, \theta_s \notin\ZZ $,
$ \frac{1}{2}(\theta_{\infty}\pm\sigma) \notin\ZZ $,
$ \frac{1}{2}(\theta_{s}\pm\theta_0\pm\sigma) \notin\ZZ $ and
that $ 0 < \Re(\sigma) < 1 $ (a distinct solution was presented for $ \sigma=0 $).
These conditions therefore strictly apply only to the generic or transcendental solutions of
the fifth Painlev\'e equation. 
For generic parameter values the terms given in (\ref{PV_sExp_0}) uniquely 
specify all the subsequent terms in the convergent Puisuex-type expansion for 
$ \zeta(s) $ about $ s=0 $
\begin{equation}
   \zeta(s) = 
   \sum^{\infty}_{j=0}\sum_{|k|\leq j}c_{j,k}s^{j+k\sigma} ,
\label{PV_puisuex}
\end{equation}
i.e. with any two of $ c_{1,0},c_{1,1} $ or $ c_{1,-1} $ given. 

To relate this to $ \tilde{E}_{2,N} $, we see from (\ref{LUE_Vtau}) and 
(\ref{LUE_theta}) that we require $ \sigma^2=(a+\mu)^2 $ and thus we can choose
\begin{equation}
  \sigma=a+\mu .
\label{LUE_const}
\end{equation}
This relation, $ \sigma=\theta_0+\theta_{s} $, is a violation of one of the strict 
conditions given above and is in fact a sufficient condition for a classical solution,
along with the necessary condition $ \theta_0+\theta_s+\theta_{\infty}=-2N\in\ZZ $,
which is the type of solution that we are dealing with here. 
However we conjecture that Jimbo's conditions
can be relaxed to accommodate such solutions and the corresponding formulae 
(or limiting forms if necessary) still hold. 
With this choice of $ \sigma $ the coefficient of $ s^{1-\sigma} $ in (\ref{PV_sExp_0})
contains a factor of
\begin{equation*}
  \frac{1}{\Gamma(\frac{\displaystyle\theta_{\infty}+\sigma}{\displaystyle 2})}
  = \frac{1}{\Gamma(-N)}
\end{equation*}
and thus vanishes. Simplifying the other terms gives
\begin{multline*}
  \tau_V(s) \sim Cs^{-N^2-N(a+\mu)}
       \Bigg\{ 1 + \frac{(2N+a+\mu)a}{2(a+\mu)}s \\
   + u \frac{\sin\pi\mu}{\sin\pi(a+\mu)}
       \frac{\Gamma(a+1)\Gamma(\mu+1)\Gamma(N+1+a+\mu)}
            {\Gamma^2(2+a+\mu)\Gamma(1+a+\mu)\Gamma(N)}
      s^{1+a+\mu} \Bigg\} .
\end{multline*}
Substituting in (\ref{LUE_Vtau}) we see that this is in precise agreement with 
(\ref{LUE_Eexp}) provided we choose
\begin{equation}
  u\frac{\sin\pi\mu}{\sin\pi(a+\mu)} = (1-\xi)e^{\pi i\mu}-\frac{\sin\pi a}{\sin\pi(a+\mu)}
\label{LUE_u}
\end{equation}

\section{The hard edge limit}\label{S4}
\setcounter{equation}{0}
The hard edge limit is defined by (\ref{HE_Edefn}). However, only in the cases
$ \mu=0 $, $ \mu=2 $ do we know how to prove its existence for general $ \xi $
(in the case $ \mu=0 $ $ \tilde{E}_{2,N} $ can be written as a Fredholm determinant,
while the case $ \mu=2 $ is related to this via differentiation). However a log-gas
viewpoint (\cite{rmt_Fo}) indicates that the limit will be well defined, and
moreover we expect that it can be taken term-by-term in the small $ s $ expansion
of $ \tilde{E}_{2,N} $. In this section we will show that taking the hard edge 
limit of the small $ s $ expansion (\ref{LUE_Eexp}) give rise to an initial 
condition for the solution of (\ref{PIII_sigma}) consistent with that allowed
by Jimbo's theory of the small $ s $ expansion of the Painlev\'e \IIId equation. From a 
practical perspective this specifies $ \tilde{E}^{\rm hard}_{2} $ for general
values of the parameters according to (\ref{HE_E}), while from a theoretical
viewpoint it lends weight to the belief that (\ref{HE_E}) is indeed the correct
limiting evaluation for general values of the parameters.

Under the assumption that the hard edge limit can be taken term-by-term in the
small $ s $ expansion of Proposition (\ref{LUE_exp}) is immediate.

\begin{corollary}
For $ \Re(\mu)>-1 $, $ \Re(a)>-1 $ and $ \mu+a\notin\ZZ_{\geq 0} $ 
we have
\begin{multline}
  \tilde{E}^{\rm hard}_{2}(s;a,\mu;\xi) = 1-\frac{\mu}{4(a+\mu)}s+{\rm O}(s^2)
  \\
  +\frac{\Gamma(\mu+1)\Gamma(a+1)}{\Gamma^2(\mu+a+2)\Gamma(\mu+a+1)}
   \left( (1-\xi)e^{\pi i\mu}-\frac{\sin\pi a}{\sin\pi(\mu+a)} \right) 
   \left(\frac{s}{4}\right)^{\mu+a+1}\left\{1+{\rm O}(s)\right\}
  \\ +{\rm O}(s^{2(\mu+a+1)}) ,
\label{HE_Eexp}
\end{multline}
and consequently the $ \sigma $-function $ \sigma_{\rm III'}(s) $ in (\ref{HE_E.2})
has the small $ s $ expansion
\begin{multline}
  \sigma_{\rm III'}(s) = -\frac{\mu(\mu+a)}{2}+\frac{\mu}{4(\mu+a)}s+{\rm O}(s^2)
  \\
  -\frac{\Gamma(\mu+1)\Gamma(a+1)}{\Gamma(\mu+a+2)\Gamma^2(\mu+a+1)}
   \left( (1-\xi)e^{\pi i\mu}-\frac{\sin\pi a}{\sin\pi(\mu+a)} \right) 
   \left(\frac{s}{4}\right)^{\mu+a+1}\left\{1+{\rm O}(s)\right\}
  \\ +{\rm O}(s^{2(\mu+a+1)}) .
\label{HE_Sexp}
\end{multline}
\end{corollary}

Some examples of exceptional cases not covered by the preceding corollary
are the following. They are obtained by taking the hard edge limit of
(\ref{LUE_Eexp_pole0}) and (\ref{LUE_Eexp_pole1}). 
\begin{corollary}
For $ \Re(\mu)>-1 $, $ \Re(a)>-1 $ and $ \mu+a=0 $ we have 
\begin{multline}
\tilde{E}^{\rm hard}_{2}(s;a,\mu=-a;\xi) = 1
  \\
  +\Big\{
   -1+\frac{\pi a}{\sin\pi a}e^{-\pi ia}(1-\xi)
     +a[2\psi(2)+\psi(1)-\psi(1-a)-\log(s/4)] \Big\}\frac{s}{4}
  \\ +{\rm o}(s) ,
\label{HE_Eexp0}
\end{multline}
whilst for $ \mu+a=1 $ we have
\begin{multline}
\tilde{E}^{\rm hard}_{2}(s;a,\mu=1-a;\xi) = 1+(a-1)\frac{s}{4}
  \\
  +\frac{a(a-1)}{4} \Big\{
   \frac{\pi}{\sin\pi a}e^{-\pi ia}(1-\xi)
     +2\psi(3)+\psi(2)-\psi(2-a)-\log(s/4) \Big\}\left(\frac{s}{4}\right)^2
  \\ +{\rm o}(s^2) .
\label{HE_Eexp1}
\end{multline}
\end{corollary}

To compare these results to the small independent variable expansions given by 
Jimbo in the theory of \IIId, we must first undertake some preliminary
calculations as the equation (\ref{PIII_sigma}) is not directly studied in 
\cite{Ji_1982}. Rather the equation studied is
\begin{equation}
  (t\zeta'')^2 = 4\zeta'(\zeta' -1)(\zeta - t\zeta')
  + \left( \frac{v_1+v_2}{2}-v_1\zeta'\right)^2 ,
\label{PIII_zeta}
\end{equation}
where we have identified $ \theta_0=-v_2 $, $ \theta_{\infty}=-v_1 $ 
($ \theta_0, \theta_{\infty} $ are the parameters appearing in \cite{Ji_1982}).
In terms of $ \zeta(t) $ the $ \tau $-function $ \tau_{\rm III'}(t) $ is specified
by the requirement that
\begin{equation}
  \zeta(t) = t\frac{d}{dt}\log\tau_{\rm III'}(t)+\frac{v_2^2-v_1^2}{4}+t ,
\label{PIII_zeta_tau}
\end{equation}
and it is the small $ t $ expansion of $ \tau_{\rm III'}(t) $ presented in 
\cite{Ji_1982}. Comparison of (\ref{PIII_zeta}) and (\ref{PIII_sigma}) shows that
\begin{equation}
  \zeta(t) = -\sigma_{\rm III'}(s)+\frac{v_1(v_2-v_1)}{4}+\frac{s}{4} ,
  \qquad t = \frac{s}{4} .
\label{HE_zeta}
\end{equation}
Recalling (\ref{HE_E.2}), (\ref{HE_E}), (\ref{HE_zeta}) and (\ref{PIII_zeta_tau}) 
we see
\begin{equation}
 \tilde{E}^{\rm hard}_{2}(s;a,\mu;\xi) = t^{(v_2^2-v_1^2)/4}\tau_{\rm III'}(t) .
\label{HE_tau}
\end{equation}

In \cite{Ji_1982} the most general small $ t $ expansion of $ \tau_{\rm III'}(t) $
as permitted by (\ref{PIII_zeta}) is presented. It reads
\begin{multline}
  \tau_{\rm III'}(t) = Ct^{(\sigma^2-v^2_2)/4}
       \Bigg\{ 1 + \frac{v_1v_2-\sigma^2}{2\sigma^2}t
   \\
   - u
      \frac{\Gamma^2(-\sigma)}{\Gamma^2(2+\sigma)}
      \frac{\Gamma(1+\frac{\displaystyle v_2+\sigma}{\displaystyle 2})
            \Gamma(1+\frac{\displaystyle v_1+\sigma}{\displaystyle 2})}
           {\Gamma(\frac{\displaystyle v_2-\sigma}{\displaystyle 2})
            \Gamma(\frac{\displaystyle v_1-\sigma}{\displaystyle 2})}
      t^{1+\sigma}
   \\
   - \frac{1}{u}
      \frac{\Gamma^2(\sigma)}{\Gamma^2(2-\sigma)}
      \frac{\Gamma(1+\frac{\displaystyle v_2-\sigma}{\displaystyle 2})
            \Gamma(1+\frac{\displaystyle v_1-\sigma}{\displaystyle 2})}
           {\Gamma(\frac{\displaystyle v_2+\sigma}{\displaystyle 2})
            \Gamma(\frac{\displaystyle v_1+\sigma}{\displaystyle 2})}
      t^{1-\sigma}
   \\
   + {\rm O}(|t|^{2(1-\Re(\sigma))}) \Bigg\} ,
\label{PIII_tExp_0}
\end{multline}
where as in (\ref{PV_sExp_0}) $ C $ is a normalisation, while $ u $ and $ \sigma $
are arbitrary parameters. This result was established under the assumptions that
$ \frac{1}{2}(v_1\pm\sigma) \notin\ZZ $ and $ \frac{1}{2}(v_2\pm\sigma) \notin\ZZ $
along with $ 0 < \Re\sigma < 1 $ (for $ \sigma=0 $ a distinct solution is given).

To see that this structure is consistent with (\ref{HE_Eexp}) and (\ref{HE_tau}),
recalling (\ref{HE_IIIparam}) we see that for the right hand side of 
(\ref{HE_tau}) to tend to $ 1 $ as $ t $ tends to zero we must have $ C=1 $ and
$ \sigma=\pm v_1 $. Again this is a violation of first condition given above but we
conjecture that the formulae have meaning under the following limiting procedure and
are correct. Choosing the positive sign for definiteness, and then writing 
\begin{equation*}
   \frac{u}{\Gamma(\frac{\displaystyle v_1-\sigma}{\displaystyle 2})}
  = \frac{u(v_1-\sigma)}{2\Gamma(1+\frac{\displaystyle v_1-\sigma}{\displaystyle 2})}
\end{equation*}
we see that requiring 
\begin{equation*}
  \frac{u}{2}(v_1-\sigma) \to \tilde{u}\frac{\sin\pi v_1}{\pi}
  \quad {\rm as} \quad \sigma \to v_1 ,
\end{equation*}
(\ref{PIII_tExp_0}) reads
\begin{multline}
  \tau_{\rm III'}(t) \sim t^{(v_1^2-v^2_2)/4}
       \Bigg\{ 1 + \frac{v_1v_2-v_1^2}{2v_1^2}t
   \\
   + \tilde{u}
      \frac{\sin\pi(v_1-v_2)/2}{\sin\pi v_1}
      \frac{\Gamma(1-\frac{\displaystyle v_2-v_1}{\displaystyle 2})
            \Gamma(1+\frac{\displaystyle v_2+v_1}{\displaystyle 2})}
           {\Gamma^2(2+v_1)\Gamma(1+v_1)}
      \left( \frac{t}{4}\right)^{1+v_1} \Bigg\} .
\label{PIII_tExp_0.2}
\end{multline}
Recalling again (\ref{HE_IIIparam}) and (\ref{HE_tau}) we see that this agrees
with (\ref{HE_Eexp}) provided
\begin{equation}
  \tilde{u}\frac{\sin\pi\mu}{\sin\pi(a+\mu)}
  = (1-\xi)e^{\pi i\mu}-\frac{\sin\pi a}{\sin\pi(a+\mu)} ,
\label{HE_u}
\end{equation}
(cf. (\ref{LUE_u})).

\section{Application}\label{S5}
\setcounter{equation}{0}
In a recent work relating to the application of random matrix theory to the
study of moments of the derivative of the Riemann zeta-function, Conrey,
Rubinstein and Snaith \cite{CRS_2005} obtained two asymptotic expressions
associated with the derivative of characteristic polynomials for random unitary
matrices. With $ U $ a Haar distributed element of the unitary group $ U(N) $,
and $ e^{i\theta_1},\ldots,e^{i\theta_N} $ its eigenvalues, let
\begin{equation}
  \Lambda_A(s) = \prod^N_{j=1}(1-se^{-i\theta_j}) ,
\label{UcharP}
\end{equation}
and
\begin{equation}
  {\mathcal Z}_A(s) = e^{-\pi iN/2}e^{i\sum^N_{n=1}\theta_n/2}s^{-N/2}\Lambda_A(s) ,
\label{UcharZ}
\end{equation}
(note that $ {\mathcal Z}_A(e^{i\theta}) $ is real for $ \theta $ real). In terms
of this notation, the two results from \cite{CRS_2005} are
\begin{equation}
  \langle |\Lambda'_A(1)|^{2k} \rangle_{A\in U(N)} \mathop{\sim}\limits_{N \to \infty}
  b_k N^{k^2+2k} ,
\label{LcharM.1}
\end{equation}
where
\begin{multline}
 b_k = (-1)^{k(k+1)/2} \sum^{k}_{h=0}{{k}\choose{h}}(k+h)!
 \\
 \times
 [x^{k+h}]\left( e^{-x}x^{-k^2/2}
          \det[I_{\alpha+\beta-1}(2\sqrt{x})]_{\alpha,\beta=1,\ldots,k}\right) ,
\label{LcharM.1aux}
\end{multline}
and
\begin{equation}
  \langle |{\mathcal Z}'_A(1)|^{2k} \rangle_{A\in U(N)} \mathop{\sim}\limits_{N \to \infty}
  b'_k N^{k^2+2k} ,
\label{ZcharM.1}
\end{equation}
where
\begin{equation}
 b'_k = (-1)^{k(k+1)/2} (2k)!
 [x^{2k}]\left( e^{-x/2}x^{-k^2/2}
         \det[I_{\alpha+\beta-1}(2\sqrt{x})]_{\alpha,\beta=1,\ldots,k}\right) .
\label{ZcharM.1aux}
\end{equation}
In (\ref{LcharM.1aux}) and (\ref{ZcharM.1aux}) the notation $ [x^p]f(x) $ denotes 
the coefficient of $ x^p $ in $ f(x) $.

The relevance of these formulae to the present study is that the determinant 
therein can be identified in terms of $ \tilde{E}^{\rm hard}_{2} $. Thus, we have
shown in a previous study \cite{FW_2002a} that for $ a \in \ZZ_{\geq 0} $
\begin{equation}
 \tilde{E}^{\rm hard}_{2}(s;a,\mu;\xi=1)
  = A(a,\mu)\left(\frac{2}{\sqrt{s}}\right)^{a\mu}e^{-s/4}
                           \det[I_{\mu+\alpha-\beta}(\sqrt{s})]_{\alpha,\beta=1,\ldots,a} .
\label{HE_bessel}
\end{equation}
where
\begin{equation}
  A(a,\mu) = a!\prod^a_{j=1}\frac{(j+\mu-1)!}{j!} .
\label{.1}
\end{equation}
Interchanging row $ \beta $ by row $ a-\beta+1 $ ($ \beta=1,\ldots,a $ in order) 
we see from this that
\begin{equation}
\begin{split}
  b_k  & = \frac{(-1)^k}{A(k,k)}\sum^k_{h=0} {{k}\choose{h}}(k+h)!
          [x^{k+h}]\tilde{E}^{\rm hard}_{2}(4x;k,k;\xi=1)
  \\
  b'_k & = \frac{(-1)^k}{A(k,k)}(2k)!
          [x^{2k}]\left(  e^{x/2}\tilde{E}^{\rm hard}_{2}(4x;k,k;\xi=1) \right)
\end{split}
\label{:a}
\end{equation}
Note that the Painlev\'e \IIId parameters appearing in this solution are 
$ \mu=a=k\in\NN $ and $ \mu+a=2k\in2\NN $ and thus we are dealing with the exceptional 
case of indeterminacy referred to in Section 2. However as was noted there the 
generic formulae still hold with to the modifications discussed and in particular 
the $\sigma$-function has a small argument expansion of a purely analytic form. 

From (\cite{FW_2003b}) it is known that the determinants in (\ref{LcharM.1aux})
and (\ref{ZcharM.1aux}) can also be expressed as a particular generalised 
hypergeometric function. Such an observation implies, for instance, that
\begin{equation}
 x^{-k^2/2}\det[I_{\alpha+\beta-1}(2\sqrt{x})]_{\alpha,\beta=1,\ldots,k} 
   = \prod^k_{j=1}\frac{j!}{\Gamma(j+k)}
     {{}^{\vphantom{(1)}}_0}F^{(1)}_1(;2k;x_1,\ldots,x_k)|_{x_j=x} ,
\label{.b}
\end{equation} 
where $ {{}^{\vphantom{(1)}}_0}F^{(1)}_1(;c;x_1,\ldots,x_k) $ has a series
development about $ x_1,\ldots,x_k=0 $ with an explicitly given coefficient for an
arbitrary term.
However this is not a practical or efficient way to compute the coefficients 
required in (\ref{LcharM.1aux}) or (\ref{ZcharM.1aux}) for moderate or large $ k $
as it involves the hook lengths of Young diagrams associated with the partitions
of $ k $. 

According to (\ref{HE_E}), (\ref{PIII_sigma}) and (\ref{HE_Sexp})
\begin{equation}
 \tilde{E}^{\rm hard}_{2}(4x;k,k;\xi=1) 
 = \exp\left( -\int^{4x}_0\frac{ds}{s}\;(\sigma_{\rm III'}(s)+k^2)\right) ,
\label{Dzeta_tau}
\end{equation}
where $ \sigma_{\rm III'}(s) $ satisfies the particular $ \sigma $-Painlev\'e 
\IIId equation
\begin{equation}
  (s\sigma''_{{\rm III}'})^2 
  + \sigma'_{{\rm III}'}(4\sigma'_{{\rm III}'}-1)(\sigma_{{\rm III}'}-s\sigma'_{{\rm III}'})
  - \frac{k^2}{16} = 0 ,
\label{Dzeta_PIIIsigma}
\end{equation}
subject to the boundary condition
\begin{equation}
  \sigma_{\rm III'}(s) \mathop{\sim}\limits_{s \to 0} -k^2+\frac{s}{8}+{\rm O}(s^2),
  \quad k\in\NN .
\label{Dzeta_PIIIBC}
\end{equation}
Substituting 
\begin{equation}
  \sigma_{\rm III'}(s) = \eta(s)+\frac{s}{8} ,
\label{Dzeta_xfm}
\end{equation}
(\ref{Dzeta_PIIIsigma}) reads
\begin{equation}
  (s\eta'')^2 + 4((\eta')^2-\frac{1}{64})(\eta-s\eta') - \frac{k^2}{4^2} = 0 .
\label{Dzeta_PIII.2}
\end{equation}
We see immediately that $ \eta(s) $ can be expanded in an even function of $ s $ about
$ s=0 $,
\begin{equation}
  \eta(s) = \sum^{\infty}_{n=0}c_ns^{2n}, \qquad c_0=-k^2,
  \quad k\in\NN .
\label{Dzeta_Exp}
\end{equation}
Moreover the coefficients can be computed by a recurrence relation.

\begin{proposition}\label{Dzeta_recur}
Substituting (\ref{Dzeta_Exp}) in (\ref{Dzeta_PIII.2}) shows
\begin{equation}
  c_1 = \frac{1}{64(4k^2-1)} ,
\label{Dzeta_c1}
\end{equation}
while for $ p \geq 2 $
\begin{multline}
 c_p = \frac{1}{2c_1p(2p-1)+(2p-1)/64-8pk^2c_1}
 \\
 \times\Bigg( 4k^2\sum^{p-2}_{l=1}(l+1)(p-l)c_{l+1}c_{p-l}
 \\
             -\sum^{p-2}_{l=1}(l+1)(p-l)(2l+1)(2p-2l+1)c_{l+1}c_{p-l}
 \\
             -\sum^{p-1}_{l=1}(1-2l)c_{l}A_{p-l-1} \Bigg) ,
\label{Dzeta_cp}
\end{multline}
where
\begin{equation}
  A_q = \sum^{q}_{l=0}(l+1)(q-l+1)c_{l+1}c_{q-l+1} .
\label{Dzeta_aux}
\end{equation}
\end{proposition}
\begin{proof}
With $ h_l := (l+1)(2l+1)c_{l+1} $ we see
\begin{equation}
  (s\eta'')^2 = 4\sum^{\infty}_{p=1}H_{p-1}s^{2p}, \qquad
  H_p = \sum^{p}_{l=0}h_lh_{p-l} ,
\label{Dzeta_aux1}
\end{equation}
and similarly with $ a_l := (l+1)c_{l+1} $ we have
\begin{equation*}
  (s\eta')^2 = 4s^2\sum^{\infty}_{p=0}A_{p}s^{2p}, \qquad
  A_p = \sum^{p}_{l=0}a_la_{p-l} .
\end{equation*}
It follows from this latter result that
\begin{equation}
  \left( (\eta')^2-\frac{1}{64} \right)(\eta-s\eta')  = \sum^{\infty}_{p=0}G_{p}s^{2p} ,
\label{Dzeta_aux2}
\end{equation}
where
\begin{equation*}
  G_p = \sum^{p}_{l=0}(1-2l)c_lb_{p-l}, \quad
  b_0 = -\frac{1}{64}, \quad b_p = 4A_{p-1} \quad (p\geq 1) .
\end{equation*}
Substituting (\ref{Dzeta_aux1}) and (\ref{Dzeta_aux2}) in (\ref{Dzeta_PIII.2})
and equating like coefficients of $ s^{2p} $ to zero shows that for $ p\geq 1 $
\begin{equation*}
  H_{p-1}+G_{p} = 0 .
\end{equation*}
This for $ p=1 $ implies (\ref{Dzeta_c1}), and for $ p>1 $ implies (\ref{Dzeta_cp}).
\end{proof}

Using Proposition \ref{Dzeta_recur} it is straightforward to calculate, via 
computer algebra, the first $ k $ coefficients in (\ref{Dzeta_Exp}) for any 
particular value of $ k $. Furthermore use of computer algebra gives the 
power series up to $ x^{2k} $ of 
\begin{equation*}
  \tilde{E}^{\rm hard}_{2}(4x;k,k;\xi=1) \quad{\rm and}\quad
  e^{x/2}\tilde{E}^{\rm hard}_{2}(4x;k,k;\xi=1) ,
\end{equation*}
according to (\ref{Dzeta_tau}). From these power series the formulae 
(\ref{:a}) are used to compute $ b_k $ and $ b'_k $. In \cite{CRS_2005} the
first 15 values of both $ b_k $ and $ b'_k $ were tabulated. This can be rapidly
extended using the present method. However the resulting rational numbers 
quickly become unwieldy to record. Let us then be content by presenting just the
16th member of the sequences,
\begin{equation*}
  b_{16} = \frac{\scriptstyle 
307
\cdot
23581
\cdot
92867
\cdot
760550281759
             }
             {\scriptstyle 
2^{272}
\cdot
3^{130}
\cdot
5^{66}
\cdot
7^{42}
\cdot
11^{24}
\cdot
13^{21}
\cdot
17^{16}
\cdot
19^{14}
\cdot
23^{10}
\cdot
29^{6}
\cdot
31^{5}
\cdot
37^{3}
\cdot
41^{2}
\cdot
43^{2}
\cdot
47
\cdot
53
\cdot
59
\cdot
61
             } ,
\end{equation*}
\begin{equation*}
  b'_{16} = \frac{\scriptstyle 
4148297603
\cdot
7623077808870586151748455369217213506671334530597
              }
              {\scriptstyle
2^{264}
\cdot
3^{133}
\cdot
5^{66}
\cdot
7^{42}
\cdot
11^{25}
\cdot
13^{21}
\cdot
17^{16}
\cdot
19^{14}
\cdot
23^{11}
\cdot
29^{7}
\cdot
31^{6}
\cdot
37^{3}
\cdot
41^{2}
\cdot
43^{2}
\cdot
47
\cdot
53
\cdot
59
\cdot
61
              } .
\end{equation*}

\section{Acknowledgements}
This work was supported by the Australian Research Council. PJF thanks M. Rubinstein
for relating the results of \cite{CRS_2005} before publication, and for the 
organisers of the Newton Institute program `Random matrix approaches in number theory'
held in the first half of 2004 for making this possible. NSW wishes to thank the 
organisers of the CRM program `Random Matrices, Random Processes and Integrable Systems'
held in Montreal 2005 for the opportunity to attend the program and in particular the 
hospitality of John Harnad during his stay.

\bibliographystyle{plain}
\bibliography{moment,random_matrices,nonlinear}
                                                                                                    
\end{document}